\documentclass[a4paper,11pt,reqno]{amsart}

\usepackage[colorlinks=true, linkcolor=blue]{hyperref}
\usepackage{cite}

\usepackage{latexsym,amsmath,amssymb,amsfonts,amsthm}%orcidlink
\usepackage{mathrsfs}
\usepackage{tikz}

\pgfarrowsdeclare{<<<}{>>>}
{
\arrowsize=0.2pt
\advance\arrowsize by .5\pgflinewidth
\pgfarrowsleftextend{-4\arrowsize-.5\pgflinewidth}
\pgfarrowsrightextend{.5\pgflinewidth}
}
{
\arrowsize=0.2pt
\advance\arrowsize by .5\pgflinewidth
\pgfpathmoveto{\pgfpointorigin}
\pgfpathlineto{\pgfpoint{-1.2mm}{-.8mm}}
\pgfusepathqstroke
\pgfpathmoveto{\pgfpointorigin}
\pgfpathlineto{\pgfpoint{-1.2mm}{.8mm}}
\pgfusepathqstroke
}%%%% DEFINIZIONE FRECCIA CARINA:FINE (per usarla, iniziare una figura con \begin{tikzpicture}[>=>>>] )

\usetikzlibrary{arrows}
\usepackage{enumerate}

\usepackage{a4wide}

\numberwithin{equation}{section}
\newtheorem{theorem}{Theorem}[section]
\newtheorem{lemma}[theorem]{Lemma}

\theoremstyle{remark}
\newtheorem{remark}[theorem]{Remark}
\theoremstyle{definition}
\newtheorem{defin}[theorem]{Definition}

\font\script=rsfs10 at 11pt
\def\H{{\mbox{\script H}\,\,}}
\def\eps{\varepsilon}
\def\L{\mathscr{L}}
\def\E{\mathcal E}
\def\P{\mathcal P}
\def\R{\mathbb R}

\def\Z{\mathcal Z}
\def\W{\mathcal W}
\def\bal{\begin{aligned}}
\def\eal{\end{aligned}}
\def\proofof#1{\begin{proof}[Proof of #1]}

\def\case#1#2{\par\noindent{\underline{\it Case~#1.}}\emph{ #2}\\}

\def\XXint#1#2#3{{\setbox0=\hbox{$#1{#2#3}{\int}$} \vcenter{\vspace{-1pt}\hbox{$#2#3$}}\kern-.5\wd0}}

\def\XXiint#1#2#3{{\setbox0=\hbox{$#1{#2#3}{\iint}$} \vcenter{\vspace{-1pt}\hbox{$#2#3$}}\kern-0.5\wd0}}

\def\comp{\subset\subset}
\def\spt{{\rm spt}}
\newcommand{\norm}[1]{\left\lVert#1\right\rVert}
\DeclareMathOperator*{\argmin}{\arg\,\min}
\DeclareMathOperator*{\argmax}{\arg\,\max}

\DeclareMathOperator*{\esslimsup}{ess\,-lim\,sup}
\DeclareMathOperator*{\essliminf}{ess\,-lim\,inf}
\DeclareMathOperator{\sign}{sgn}

\newcommand\restr[2]{\ensuremath{\left.#1\right|_{#2}}}

\pgfarrowsdeclare{<<<}{>>>}
{
\arrowsize=0.2pt
\advance\arrowsize by .5\pgflinewidth
\pgfarrowsleftextend{-4\arrowsize-.5\pgflinewidth}
\pgfarrowsrightextend{.5\pgflinewidth}
}
{
\arrowsize=0.2pt
\advance\arrowsize by .5\pgflinewidth
\pgfpathmoveto{\pgfpointorigin}
\pgfpathlineto{\pgfpoint{-1.2mm}{-.8mm}}
\pgfusepathqstroke
\pgfpathmoveto{\pgfpointorigin}
\pgfpathlineto{\pgfpoint{-1.2mm}{.8mm}}
\pgfusepathqstroke
}%%%% DEFINIZIONE FRECCIA CARINA:FINE (per usarla, iniziare una figura con \begin{tikzpicture}[>=>>>] )

\pgfarrowsdeclare{<<<<}{>>>>}
{
\arrowsize=0.2pt
\advance\arrowsize by .5\pgflinewidth
\pgfarrowsleftextend{-4\arrowsize-.5\pgflinewidth}
\pgfarrowsrightextend{.5\pgflinewidth}
}
{
\arrowsize=0.2pt
\advance\arrowsize by .5\pgflinewidth
\pgfpathmoveto{\pgfpointorigin}
\pgfpathlineto{\pgfpoint{-1.8mm}{-1.2mm}}
\pgfusepathqstroke
\pgfpathmoveto{\pgfpointorigin}
\pgfpathlineto{\pgfpoint{-1.8mm}{1.2mm}}
\pgfusepathqstroke
}%%%% DEFINIZIONE GROSSA:FINE (per usarla, iniziare una figura con \begin{tikzpicture}[>=>>>>] )

\newcounter{mt}
\def\maintheorem#1#2#3{\par \medskip \noindent {\bf Theorem~\mref{#1}}~(#2).~{\it #3}\par}
\def\mref#1{\Alph{#1}}
\def\maintheoremdeclaration#1{\stepcounter{mt}\newcounter{#1}\setcounter{#1}{\arabic{mt}}}
\maintheoremdeclaration{continuity}
\allowdisplaybreaks

\title[Continuity of critical points of non-local energies in 1D]{Continuity of critical points\\for 1-dimensional non-local energies}

\author{D. Carazzato}\address[D. Carazzato]{Faculty of Mathematics, University of Vienna, Austria}\email{davide.carazzato@univie.ac.at
}
\author{N. Fusco}\address[N. Fusco]{Dipartimento di Matematica, Universit\`a di Napoli, Italy}\email{n.fusco@unina.it}
\author{A. Pratelli}\address[A. Pratelli]{Dipartimento di Matematica, Universit\`a di Pisa, Italy}\email{aldo.pratelli@unipi.it
}

\date{\today}

\begin{document}
\begin{abstract}
In this paper we deal with the bounded critical points of a Riesz energy of attractive-repulsive type in dimension 1. Under suitable assumptions on the growth of the kernel in the origin, we are able to prove that they are continuous inside their support.
\end{abstract}
    
\maketitle

\section{Introduction}
In the recent years the mathematical community has been more and more interested about energies related to attractive-repulsive forces. To be more specific, we are interested in a functional of convolution type, that in its maximal generality can be defined in the class of probability measures $\P(\R^d)$:
\begin{equation}\label{eq:problem}
\E(\mu) = \iint g(x-y) d \mu(x) d \mu(y)\,,
\end{equation}
where $g:\R^d\to[0,+\infty]$ is the kernel describing the interaction. For instance, when $g$ is radial and it is not radially monotone, then the minimization of $\E$ is a non-trivial problem, and can lead to a large variety of behaviors. Notice that in several models a part of the total energy $\E$ is given by a monotone kernel. See for example some works about different versions of the Gamow model~\cite{KM2013,KM2014,F2M3,CFP2023,C2023}. We also mention~\cite{F2023} for a comparison between different attractive-repulsive models containing the Riesz energy. The prototypical kernel $g$ that describes an attractive-repulsive interaction is of the form
\[
g(x) = \frac{|x|^{\alpha}}{\alpha}-\frac{|x|^{\lambda}}{\lambda}\,,
\]
where $\alpha>0$ and $-d<\lambda<\alpha$, and many questions are still open even for these simple functions. This kind of energies appear in many different models, and they have been intensively studied both from the analytical and the numerical point of view. In fact, some very general existence result is available (see Section~\ref{sec:preliminary-results}), whereas only in some very special cases a precise formula for the minimizers has been found (see~\cite[Section 5]{CS2023}, \cite{F2022,FM2025}). Their shape changes even more dramatically when the kernel is not singular in the origin (see~\cite{CFP2017,DLM2022-1,DLM2022-2,CPT2024}). Instead, some regularity results valid for singular kernels are contained in~\cite{CDM2016}, where they exploit the specific singularity of the kernel, related to the fractional Laplacian, to address the problem via PDE techniques. Our starting point is~\cite{CP2025}, where, loosely speaking, the authors assume that $g$ satisfies a differential inequality and obtain an $L^{\infty}$ estimate for the minimizers of $\E$ (see Theorem~\ref{prelres} in Section~\ref{sec:preliminary-results} for the precise result).

In this article, we obtain the internal continuity of the bounded local minimizers of $\E$; more precisely, we consider measures which are not necessarily minimizers, but which have constant potential inside their support (see~\eqref{eq:potential} for the definition of potential). This is more general, since the constancy of the potential is ensured for all the local minimizers by the Euler-Lagrange conditions. 

We stress that we are interested in the \emph{inner regularity}, i.e. we will deal only with points in the interior of the support of a measure. In fact, in some cases minimizers are known to be bounded away from $0$ near boundary points of their support (see for example the discussion in~\cite{CDM2016} and the anisotropic version contained in~\cite[Theorem~1.1]{CMMRSV2020}, which are related to kernels). We also refer to~\cite{FM2025} for some explicit minimizers, that are unbounded in some cases. Some partial information about the boundary behaviour is available, only for minimizers, applying~\cite[Lemma 3.11]{CP2025}, where they exploit the complete Euler-Lagrange conditions~\eqref{eq:EL}. Our main theorem, namely Theorem~\mref{continuity}, regards a kernel $g$ that satisfies some hypotheses, stated below, and considers optimal measures which are already known to be $L^\infty$. This extra assumption is in fact very mild, since optimal measures are actually known to be bounded under assumptions which are only slightly stronger, as we will describe in Section~\ref{sec:preliminary-results}. Moreover, it is known that some hypotheses about the behaviour of $g$ near the origin are necessary to have some regularity. In particular, if the kernel $g$ is not singular, then we may also have some minimizers that are not absolutely continuous with respect to the Lebesgue measure (see for example~\cite{CFP2017, DLM2022-1,DLM2022-2}). We remark, however, that the singularity of the kernel is not sufficient to have the absolute continuity of the minimizers, as it is showed in~\cite[Theorem 1]{FM2025}. Furthermore, even if the minimizers admit a density with respect to $\L^d$, they may not be bounded, while the control in $L^{\infty}$ is necessary in our proof (see~\cite[Section 5]{CS2023} and~\cite[Theorem 2]{FM2025} where they even find explicitly the minimizers for some special kernels). Notice, however, that their explicit minimizers are always bounded when the space dimension $d$ is 1. In this sense, we observe that they are not an obstruction to our 1-dimensional study, while in principle it is possible that there exists a kernel $g$ for which the energy $\E$ admits unbounded minimizers in dimension 1.

\subsection{Setting}
In this paper, we consider a kernel $g$ satisfying the following assumptions:
\begin{enumerate}
\item $g$ belongs to $C^1\big(\R\setminus \{0\}\big)\cap L^1_{loc}(\R)$, it is symmetric, and it is decreasing in $(0,r)$ for some length-scale $r>0$;
\item $g$ is convex in $(0,r)$, $g'$ is concave in $(0,r)$, $g'\in BV_{loc}((0,+\infty))$, and there exists $\Lambda\in(1,+\infty]$ such that
\begin{equation}\label{defLambda}
\liminf_{x\to 0}\frac{|g'(x/2)|}{|g'(x)|} = \Lambda\,.
\end{equation}
\end{enumerate}

\begin{remark}
We notice that any prototypical kernel of the form
\[
g(x) = \frac{|x|^{\alpha}}{\alpha} -\frac{|x|^{\lambda}}{\lambda}\,,
\]
for $\alpha>0$ and $-1<\lambda<\min\{1,\alpha\}$ satisfies our hypotheses, where the expression $|x|^{\lambda}/\lambda$ is intended to be $\log |x|$ when $\lambda=0$.
\end{remark}
Our main result is the following theorem (we refer to Section~\ref{sec:preliminary-results} for the definition of the potential $\psi_f$).

\maintheorem{continuity}{Continuity of critical functions}{Let $g:\R\to[0,+\infty]$ be a function satisfying hypotheses \emph{(i)-(ii)}. If $f\L^1 \in\mathcal{M}(\R)$ is a finite measure with $\spt f = [a,b]$ for some $a,b\in\R$, $\norm f_{\infty} = M<+\infty$ and $\psi_f$ is constant almost everywhere in $[a,b]$, then $f$ is continuous in $(a,b)$.}

\subsection*{Structure of the paper} In Section~\ref{sec:preliminary-results} we introduce some notation and we expose some known results. We also provide the definition of essential limits in Definition~\ref{def:essential-limit}, that are the ``Lebesgue compliant'' version of the usual limits, which are necessary since we work with a density that is well defined almost everywhere.\\
In Section~\ref{sec:second-order} we introduce the tools to look at the second derivative of a potential, since our strategy is to show that the second derivative of the potential cannot be $0$ at some point $\bar{x}$ if $f$ does not admit essential limit at $\bar{x}$. This procedure needs a preliminary step, where we regularize $f$ by convolution, and this is very convenient in our problem since the potential behaves well in terms of convolutions. In Section~\ref{subsec:cancelling}, we use the hypotheses on $g$ to show some cancellation phenomena, yielding to Lemma~\ref{lemma:convex-cancellation} and Lemma~\ref{lemma:concave-cancellation}. With these results at hand, it is not hard to prove Lemma~\ref{lemma:rearrangement-comparison}, that provides the leading term in our estimates of the main theorem. However, one has to be careful: the expression that we use for the second derivative is obtained from an integration by parts, and this is not always well defined even if the density is smooth because $g'$ is not locally integrable; more precisely, our formulas work only if one considers intervals starting at critical points of the fuction. As a consequence, we will estimate the potential only in some critical points of the regularized function $f_\delta$; Section~\ref{subseq:selection-critical-points} will be devoted to the careful choice of some critical points. The final Section~\ref{sec:proof} is then devoted to the proof of Theorem~\mref{continuity}, where we combine the various estimates found before to show that the second derivative of $f_\delta$ cannot be constantly $0$ if the convolution parameter is small enough.

\subsection{Preliminary results\label{sec:preliminary-results}}

This class of problems has been intensively studied, and some known results are the starting point for our analysis. To begin with, there exists a minimizer, and this was proved for example in~\cite{SST2015,CCP2015,CP2025}, where very weak assumptions are needed. Moreover, it was recently proved that any minimizer is absolutely continuous with respect to $\L^1$ and its density is an $L^{\infty}$ function:
\begin{theorem}[\!\!{\cite[Proposition~3.4]{CP2025}}]\label{prelres}
Let $g:\R\to[0,+\infty]$ be a function satisfying \emph{(i)} and, in addition, $g$ is strictly convex in $(0,+\infty)$. Then any minimizer $\mu\in\P(\R)$ of~\eqref{eq:problem} is of class $L^{\infty}$ and $\spt \mu = [a,b]$ for some $a,b\in\R$.
%We underline that $b-a$ and $\norm{\mu}_{\infty}$ are controlled from above by a constant depending only on the kernel $g$.
\end{theorem}

We define now a fundamental tool, the so-called \emph{potential} induced by a measure, which is defined pointwise by
\begin{equation}\label{eq:potential}
    \psi_{\mu}(x) = \int_\R g(x-y)d \mu(y).
\end{equation}
Of course, this is nothing else than $g*\mu$, it is lower semicontinuous whenever so is $g$, and $\E(\mu) = \int \psi_{\mu}(x)d \mu(x)$. This function plays a crucial role since it appears in the Euler-Lagrange conditions satisfied by a minimizer $\mu$ (see for instance~\cite{BCT2018,CDM2016,CP2025}):
\begin{equation}\label{eq:EL}
\begin{cases}
\psi_{\mu} = \E(\mu) & \mu-a.e.\,,\\
\psi_{\mu}\geq \E(\mu) & in\  \R\setminus \spt \mu\,.
\end{cases}
\end{equation}
In the sequel, when working with an absolutely continuous measure $\mu=f \L^1$, we will write $\psi_f$ in place of $\psi_\mu$.

\begin{defin}[Essential directional limits]\label{def:essential-limit}
Given a function $F:\R\to\R$ and $\bar x\in\R$, we say that $l\in\R$ is the \emph{essential liminf from the left} of $F$ at $\bar x$ if for every $\eps>0$ there exists $\eta>0$ such that
\[
\begin{cases}
\L^1(\{F>l-\eps\}\cap (\bar{x}-\eta,\bar{x})) = \eta\,,\\
\L^1(\{F<l+\eps\}\cap (\bar x-\eta,\bar x))>0\,.
\end{cases}
\]
In this case, we write $l = \essliminf_{t\to \bar x^-}F(t)$. Similar definitions can be given for the essential limsup and for the limits from the right. If all of the four essential limits coincide, then we say that $F$ admits essential limit at $\bar x$.
\end{defin}

\section{Second derivative of the potential and selection of critical points\label{sec:second-order}}

Roughly speaking, our goal is to study the second derivative of $\psi_f$, and to check that it cannot be $0$ when $f$ is not continuous. In order to make this in a formally correct way, we will regularize $f$ so to work with a smooth potential. To do that, we fix a symmetric mollifier $\rho\in C^{\infty}_c(\R)$ such that $\spt \rho = [-1,1]$, $\norm{\rho}_1=1$, $\norm{\rho}_{\infty}\leq 1$, and $\rho'< 0$ in $(0,1)$. Then, we set as usual $\rho_\delta(t) = \rho(t/\delta)/\delta$, and we consider the smooth function $f_\delta = f*\rho_\delta$, which has compact support; moreover, since the potential $\psi_f$ is constant in $(a,b)$, then the potential $\psi_{f_\delta}$ is constant in $(a+\delta,b-\delta)$. The proof will come by the fact that, if $f$ is not continuous, then $\psi_{f_\delta}''$ cannot be $0$ at some critical point of $f_\delta$. The reason to consider a critical point is merely technical; in fact, we will use an alternative formula for $\psi_{f_\delta}''$ obtained integrating by parts, that can be justified only if the point where we compute that derivative is a critical point of $f_\delta$. The contradiction will then provide the thesis.\par

Let us be more precise. Let $F\in C^{\infty}_c(\R)$ be a generic smooth function, and let $x$ be a critical point for $F$. Then, we have the following expression for $\psi_ F''(x)$:
\begin{equation}\label{eq:second-derivative}
\psi_F''(x) = \int_\R F''(t)g(x-t)\,dt = -\int_{\R} F'(t) \frac{d}{dt}g(t-x)\,dt= \int_{\R}\sign(x-t)F'(t)g'(|x-t|)\,dt\,.
\end{equation}
The fact that the integrals containing first derivatives are well defined is true since $x$ is a critical point of $F$ (which implies  $|F'(t)|\lesssim |t-x|$ since $F$ is smooth), and in view of the following lemma.
\begin{lemma}\label{lemma:integrability-g'}
Let $g:\R\to[0,+\infty]$ be a function in $L^1_{loc}(\R)\cap C^1(\R\setminus \{0\})$ such that $g'\leq 0$ in an interval $(0,r)$. Then
\[
-\int_0^r g'(t)t < +\infty\,.
\]
\end{lemma}
\begin{proof}
For any $t\in(0,r)$ we write $g(r)-g(t) = \int_t^r g'(s)\,ds$. Then, since $g'$ has constant sign in $(0,r)$, we can apply Fubini Theorem and see that
\[
\int_0^r g(r)-g(t)\, dt = \int_0^r\, dt \int_t^r g'(s)\,ds = \int_0^r\,ds \int_0^s g'(s)\, dt = \int_0^r sg'(s)ds\,.
\]
This concludes the proof since $g\in L^1([0,r])$.
\end{proof}

\subsection{Cancellation lemmas\label{subsec:cancelling}}
We are going to manipulate the expression~\eqref{eq:second-derivative}, and we obtain some inequalities for the contribution due to the integral between two critical points in that expression.

\begin{lemma}\label{lemma:convex-cancellation}
Let $g:\R\to[0,+\infty]$ satisfy \emph{(i)-(ii)}. Let $\alpha,\beta\in\R$ be given, with $\alpha<\beta$ and $\beta-\alpha<r$. Let $F\in C^2([\alpha,\beta])$ be a function with $F(\alpha) = F(\beta)$, and such that $\alpha$ and $\beta$ are absolute minimum points of $F$ in $[\alpha,\beta]$. If $F'(\beta) = 0$, then for every $x\geq \beta$ with $x-\alpha < r$ we have that
\[
\int_\alpha^\beta F'(t)g'(x-t)\,dt \geq 0\,.
\]
If, instead, $F'(\alpha) = 0$, then for every $x\leq \alpha$ with $\beta-x < r$ we have that
\[
\int_\alpha^\beta -F'(t)g'(t-x)\,dt \geq 0\,.
\]
\end{lemma}
\begin{proof}
First of all, we observe that Lemma~\ref{lemma:integrability-g'} guarantees that the integral is finite. It is also easy to check that the second inequality can be deduced from the first one considering the function $G(t)=F(-t)$ defined in the interval $[-\beta,-\alpha]$, so we will just prove the first one. Moreover, it is sufficient to prove the result when $F'$ changes sign only once: the general result can be obtained approximating $F$ with functions whose derivative has a finite number of sign changes (see the proof of Lemma~\ref{lemma:rearrangement-comparison}, where this procedure is slightly more complex). Therefore, we need to prove the result when there exists $\xi\in(\alpha,\beta)$ such that $F_1 = \restr F{(\alpha,\xi)}$ is monotone increasing, $F_2 = \restr F{(\xi,\beta)}$ is monotone decreasing, and $F'\neq 0$ in $(\alpha,\beta)\setminus \{\xi\}$. With this reduction, we use the change of variables $z=F_1(t)$ and $w=F_2(t)$ to get that
\[\begin{split}
\int_{\alpha}^{\beta}F'(t)g'(x-t)\,dt &= \int_{\alpha}^{\xi}F_1'(t)g'(x-t)\,dt+\int_{\xi}^{\beta}F_2'(t)g'(x-t)\,dt\\
&= \int_{F_1(\alpha)}^{F_1(\xi)}g'(x-F_1^{-1}(z))\,dz - \int_{F_2(\beta)}^{F_2(\xi)}g'(x-F_2^{-1}(w))\,dw\\
& = \int_{F(\alpha)}^{F(\xi)}\left[g'(x-F_1^{-1}(z))-g'(x-F_2^{-1}(z))\right]\,dz.
\end{split}\]
For every $z\in[F(\alpha),F(\xi)]$ we have that $x-F_2^{-1}(z)\leq x-F_1^{-1}(z) <r$, and since $g$ is convex in $(0,r)$, then the function inside the integral is non-negative, concluding the proof.
\end{proof}

\begin{lemma}\label{lemma:concave-cancellation}
Let $g:\R\to[0,+\infty]$ satisfy \emph{(i)-(ii)}. Let $\alpha,\beta\in\R$ be given, with $\alpha<\beta$ and $\beta-\alpha<r$. If $F\in C^2([\alpha,\beta])$, $\alpha$ is an absolute minimum point of $F$ in $[\alpha,\beta]$ and $F'(\alpha)=0$, then for every $x\leq y\leq \alpha$ with $\beta-x < r$ we have that
\[
\int_{\alpha}^{\beta} F'(t)(g'(t-x)-g'(t-y))\,dt \geq 0\,.
\]
\end{lemma}
\begin{proof}
We use a cancellation principle similar to the one of Lemma~\ref{lemma:convex-cancellation}, this time relying on the concavity of $g'$. As before, by approximation it is sufficient to prove the result when there exists $\xi\in(\alpha,\beta)$ such that $F_1 = \restr F{(\alpha,\xi)}$ is monotone increasing, $F_2 = \restr F{(\xi,\beta)}$ is monotone decreasing, and $F'\neq 0$ in $(\alpha,\beta)\setminus \{\xi\}$. Using the change of varibles $z=F_1(t)$ and $w=F_2(t)$ we arrive to
\[\begin{split}
\int_{\alpha}^{\beta} F'(t)(g'(t-&x)-g'(t-y))\,dt
= \int_{F_1(\alpha)}^{F_1(\xi)} g'(F_1^{-1}(z)-x)-g'(F_1^{-1}(z)-y)\,dz\\
& \qquad -\int_{F_2(\beta)}^{F_2(\xi)}g'(F_2^{-1}(w)-x)-g'(F_2^{-1}(w)-y)\,dw\\
& \geq \int_{F(\alpha)}^{F(\xi)}g'(F_1^{-1}(z)-x)-g'(F_2^{-1}(z)-x)\,dz\\
& \qquad -\int_{F(\alpha)}^{F(\xi)}g'(F_1^{-1}(z)-y)-g'(F_2^{-1}(z)-y)\,dz.
\end{split}\]
For any $z\in[F(\alpha),F(\xi)]$ we have that $F_1^{-1}(z)-x\leq F_2^{-1}(z)-x$ and $F_1^{-1}(z)-y\leq F_2^{-1}(z)-y$, and since $g'$ is concave in $(0,r)$, with $\beta-x<r$, then
\[
g'(F_1^{-1}(z)-x)-g'(F_2^{-1}(z)-x) \geq g'(F_1^{-1}(z)-y)-g'(F_2^{-1}(z)-y)\qquad \forall z\in[F(\alpha),F(\xi)],
\]
and the inequality is proved.    
\end{proof}

\begin{lemma}\label{lemma:rearrangement-comparison}
Let $g:\R\to[0,+\infty]$ satisfy \emph{(i)-(ii)}. Let $\alpha,\beta\in\R$ with $\alpha<\beta$ and $\beta-\alpha=\gamma<r$. If $F\in C^2([\alpha,\,\beta])$ with $F'(\alpha) = F'(\beta) = 0$ and $\alpha$ and $\beta$ are respectively an absolute minimum and an absolute maximum of $F$ in that interval, then
\begin{equation}\label{eq:cancellation-inequality}
\int_{\alpha}^{\beta}F'(t) \big( |g'(t-\alpha)|+|g'(\beta-t)|\big)\,dt\geq (F(\beta)-F(\alpha))\big(|g'(\gamma)|+|g'(\gamma/2)|\big)\,.
\end{equation}
\end{lemma}
\begin{proof}
First of all we observe that, since $F'=0$ in $\alpha,\,\beta$ and $F'$ is Lipschitz, then the integrals are finite thanks to Lemma~\ref{lemma:integrability-g'}. We can easily check that~(\ref{eq:cancellation-inequality}) is valid under the additional assumption that $F$ is increasing. Indeed, keeping in mind that $g$ is convex and $g'\leq 0$ in $(0,\gamma)$, we have that
\begin{align*}
\int_{\alpha}^{\beta}F'(t)|g'(t-\alpha)|\,dt &\geq \int_{\alpha}^{\frac{\alpha+\beta}2}F'(t)|g'(\gamma/2)|\,dt+\int_{\frac{\alpha+\beta}2}^{\beta}F'(t)|g'(\gamma)|\,dt\,,\\
\int_{\alpha}^{\beta}F'(t)|g'(\beta-t)|\,dt&\geq \int_{\alpha}^{\frac{\alpha+\beta}2}F'(t)|g'(\gamma)|\,dt+\int_{\frac{\alpha+\beta}2}^{\beta}F'(t)|g'(\gamma/2)|\,dt\,,
\end{align*}
and adding these two inequalities we get the thesis. Notice that this proof does not require that $F\in C^2$; indeed, it is sufficient that $F$ is a Lipschitz function such that
\begin{equation}\label{neededass}
|F'(\alpha+s)|+|F'(\beta-s)| \leq Cs\,, \qquad \hbox{for some $C>0$ and a.e. $\bal 0\leq s \leq \frac{\beta-\alpha}2\eal$}\,.
\end{equation}
Let us now consider the general case when $F$ is not necessarily monotone. Assume first that, as in Figure~\ref{fig:monotone-replacement}, there are $\alpha\leq \alpha' <\gamma'<\beta'\leq (\alpha+\beta)/2$ such that $F(\alpha')=F(\beta')$, and $F$ is increasing in $[\alpha',\gamma']$ and decreasing in $[\gamma',\beta']$, and call $\widetilde F$ the function, depicted in red in the figure, which equals $F$ outside the interval $[\alpha',\beta']$, and which is constantly equal to $F(\alpha')=F(\beta')$ in  $[\alpha',\beta']$.
\begin{figure}[thbp]
\begin{tikzpicture}[>=>>>] 
\draw[->] (-.5,0)--(8.5,0);
\draw[->] (0,-.5)--(0,4);
\fill (0.5,0) circle (2pt);
\draw (0.5,0) node[anchor=north] {$\alpha$};
\fill (4,0) circle (2pt);
\draw (4,0) node[anchor=north] {$\frac{\alpha+\beta}2$};
\fill (7.5,0) circle (2pt);
\draw (7.5,0) node[anchor=north] {$\beta$};
\draw[line width=.8pt] (.5,0.5) .. controls (1.5,0.5) and (2,3.5) .. (2.5,2.5) .. controls (3,1.5) and (4,1.5) .. (4.5,2.5) .. controls (5,3.5) and (6,4) .. (7.5,4);
\draw[dashed] (1.76,0)--(1.76,2);
\draw[dashed] (2.88,0)--(2.88,2);
\draw[dashed] (2.28,2.7)--(2.28,0);
\draw[red] (1.76,2)--(2.88,2);
\draw[red] (2.6,2) node[anchor=north] {$\widetilde F$};
\draw (1.5,3) node[anchor=north] {$F$};
\fill (1.76,0) circle (2pt);
\draw (1.76,0) node[anchor=north] {$\alpha'$};
\fill (2.88,0) circle (2pt);
\draw (2.88,0) node[anchor=north] {$\beta'$};
\fill (2.28,0) circle (2pt);
\draw (2.28,0) node[anchor=north] {$\gamma'$};
\end{tikzpicture}
\caption{A possible function $F$ and its replacement $\widetilde F$ in red for~(\ref{tre}).}
\label{fig:monotone-replacement}
\end{figure}
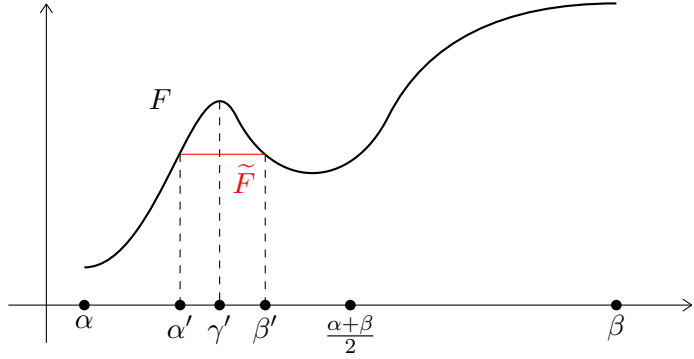
Arguing as in the proofs of the previous lemmas with a change of variable, and calling this time $F_1$ and $F_2$ the restrictions of $F$ to $[\alpha',\gamma']$ and to $[\gamma',\beta']$ respectively, since $g'\leq 0$ in $(0,r)$ we have
\begin{equation}\label{uno}\begin{split}
\int_{\alpha'}^{\beta'} F'(t) |g'(t-\alpha)|\, dt &= 
\int_{\alpha'}^{\gamma'} F'(t) |g'(t-\alpha)|\, dt+\int_{\gamma'}^{\beta'} F'(t) |g'(t-\alpha)|\, dt\\
&=\int_{F(\alpha')}^{F(\gamma')} |g'(F_1^{-1}(s)-\alpha)|\,ds-\int_{F(\alpha')}^{F(\gamma')} |g'(F_2^{-1}(s)-\alpha)|\,ds\,.
\end{split}\end{equation}
In the very same way,
\begin{equation}\label{due}
\int_{\alpha'}^{\beta'} F'(t) |g'(\beta-t)|\, dt = \int_{F(\alpha')}^{F(\gamma')} |g'(\beta-F_1^{-1}(s))|\,ds-\int_{F(\alpha')}^{F(\gamma')} |g'(\beta-F_2^{-1}(s))|\,ds\,.
\end{equation}
Let us now fix any $F(\alpha')<s<F(\gamma')$, and call for brevity $z=F_1^{-1}(s)$ and $w=F_2^{-1}(s)$. Then, since $z\leq w\leq (\alpha+\beta)/2$, we have
\[
z-\alpha \leq w-\alpha\leq \beta-w \leq \beta-z\,.
\]
As a consequence, using again the convexity of $|g'|$ in $(0,r)$, we have
\[
|g'(z-\alpha)| -|g'(w-\alpha)| + |g'(\beta-z)|-|g'(\beta-w)|\geq 0\,.
\]
Inserting this inequality in~(\ref{uno}) and~(\ref{due}), and keeping in mind that $\widetilde F\equiv F$ outside of the interval $[\alpha',\beta']$, while $\widetilde F'\equiv 0$ in $[\alpha',\beta']$, we obtain that
\begin{equation}\label{tre}
\int_{\alpha}^{\beta}F'(t)\big( |g'(t-\alpha)+|g'(\beta-t)|\big)\, dt \geq \int_{\alpha}^{\beta}\widetilde F'(t)\big( |g'(t-\alpha)+|g'(\beta-t)|\big)\, dt\,.
\end{equation}
A simple approximation argument implies that~(\ref{tre}) is valid also with $\widehat F$ in place of $\widetilde F$, where
\[
\widehat F(t) = \left\{\begin{array}{cc}
\bal \inf \bigg\{F(s),\, t\leq s \leq \frac{\alpha+\beta} 2 \bigg\}\eal & \hbox{if $\bal t\leq \frac{\alpha+\beta}2\eal$}\,,\\ 
F(t) &\hbox{if $\bal t\geq \frac{\alpha+\beta}2\eal$}\,.
\end{array}\right.
\]
And then, repeating the very same argument in the second half of the interval $[\alpha,\beta]$, we get the validity of~(\ref{tre}) also for the function $F^*$ defined as
\[
F^*(t) =  \left\{\begin{array}{cc}
\bal \inf \bigg\{F(s),\, t\leq s \leq \frac{\alpha+\beta} 2 \bigg\}\eal & \hbox{if $\bal t\leq \frac{\alpha+\beta}2\eal$}\,,\\[10pt]
\bal \sup \bigg\{F(s),\, \frac{\alpha+\beta} 2\leq s \leq  t \bigg\}\eal & \hbox{if $\bal t\leq \frac{\alpha+\beta}2\eal$}\,.
\end{array}\right.
\]
We are now ready to conclude. Indeed, the validity of~(\ref{tre}) with $F^*$ in place of $\widetilde F$ obiously implies the estimate~(\ref{eq:cancellation-inequality}) for $F$ as soon as it is proved for $F^*$. And in turn, this is true because by construction $F^*$ is increasing and satisfies~(\ref{neededass}), and we have already noticed the the inequality is true in this case.
\end{proof}
\begin{remark}\label{rem:change-sign}
An obvious consequence of~(\ref{eq:cancellation-inequality}) is that at least one between $\int_\alpha^\beta F'(t) |g'(t-\alpha)|\,dt$ and $\int_\alpha^\beta F'(t) |g'(\beta-t)|\,dt$ is greater than
\[
\frac{F(\beta)-F(\alpha)}2\, \Big(|g'(\gamma)|+|g'(\gamma/2)|\Big)\,.
\]
This can be rewritten by saying that there is a point $p\in \{\alpha,\,\beta\}$ such that
\begin{equation}\label{basecase}
\pm \int_\alpha^\beta F'(t) \, \frac d{dt}\, g(|t-p|)\,dt \geq \frac{|F(\beta)-F(\alpha)|}2\, \Big(|g'(\gamma)|+|g'(\gamma/2)|\Big)\,,
\end{equation}
where the sign is $+$ if $p$ is the maximum point, and $-$ if $p$ is the minimum point. Suppose now that $F$ is a function satisfying all the assumptions of Lemma~\ref{lemma:rearrangement-comparison} except for the fact that $\alpha$ is a maximum point and $\beta$ a minimum. In this case, the estimate~(\ref{eq:cancellation-inequality}) becomes
\[
\int_{\alpha}^{\beta}F'(t) \big( |g'(t-\alpha)|+|g'(\beta-t)|\big)\,dt\leq (F(\beta)-F(\alpha))\big(|g'(\gamma)|+|g'(\gamma/2)|\big)\,.
\]
Then, also in this case we have a point $p\in\{\alpha,\,\beta\}$ such that~(\ref{basecase}) holds, and again the sign is $+$ if $p$ is the maximum point and $-$ if $p$ is the minimum point. Summarizing, the validity of~(\ref{basecase}) with a point $p\in \{\alpha,\,\beta\}$ with sign $+$ (resp.,\, $-$) if $p$ is the maximum (resp., the minimum) point is true as soon as $F\in C^2([\alpha,\,\beta])$ is a function with $F'(\alpha)=F'(\beta)=0$ and such that one among $\alpha$ and $\beta$ is an absolute maximum of $F$ in $[\alpha,\,\beta]$, and the other one an absolute minimum.
\end{remark}

\subsection{Selection of critical points\label{subseq:selection-critical-points}}

This short, technical section is devoted to fix some notation and to select some critical points of a given function $f:(a,b)\to\R$. We concentrate ourselves on a point in $(a,b)$, that we assume to be $0$ for simplicity. We will write
\begin{align*}
l_L^-= \essliminf_{t\to0^-} f(t)\,, && l_L^+ = \esslimsup_{t\to0^-}f(t)\,,\\
l_R^- = \essliminf_{t\to0^+} f(t)\,, && l_R^+ = \esslimsup_{t\to0^+}f(t)\,,
\end{align*}
as well as
\begin{align*}
h_L = l_L^+-l_L^-\,, && h_R = l_R^+-l_R^-\,.
\end{align*}
For later use, we define now the constant
\[
\overline \Lambda = \min \bigg\{2 ,\, \frac{1+\Lambda}2\bigg\} >1\,,
\]
where $\Lambda$ is the constant appearing in~(\ref{defLambda}). The rough idea of the proof is that if $f$ is not continuous, then it has to oscillate near some point, and then also the regularised function $f_\delta$ must oscillate if $\delta$ is small enough. We will then find a contradiction working on local maxima and minima of the smooth function $f_\delta$. We fix now some of these critical points, and it is enough to consider two very specific cases.

\case{I}{If $f(x)=f(-x)$ and $h_L>0$.}
The first case that we consider is when the left jump $h_L$ is strictly positive, and $f$ is symmetric. In this case, we first fix the parameters $\eps>0$ and $0<\eta<r/4$ so that
\begin{align}\label{choiceepsetaI}
\eps < \frac{(\overline \Lambda-1)h_L}{10\overline \Lambda+6} < \frac{h_L}4\,, &&
\L^1\Big(\big\{l_L^--\eps<f<l_L^++\eps\big\}\cap \big(-2\eta,0\big)\Big) = 2\eta\,.
\end{align}
Notice that the first requirement makes sense since $\overline\Lambda>1$, and the second one is satisfied by any $\eta$ small enough by definition of the essential liminf. In addition, we can require that
\begin{equation}\label{goodLambda}
\big|g'(x/2)\big| > \overline\Lambda\big| g'(x)\big| \qquad \forall\, 0<x< 2\eta\,.
\end{equation}
This is also true as soon as $\eta$ is small enough by the definition~(\ref{defLambda}) of $\Lambda$ and by the fact that $\overline\Lambda < \Lambda$.

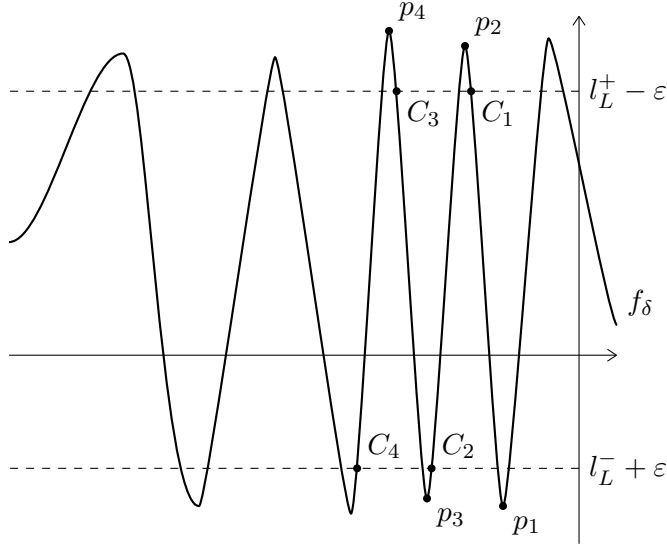
\begin{figure}[thbp]
\begin{tikzpicture}[>=>>>] 
\draw[->] (-5.5,-1)--(2.5,-1);
\draw[->] (2,-3.5)--(2,3.5);
\draw[dashed] (2,2.5)--(-5.5,2.5);
\draw (2,2.5) node[anchor=west] {$l_L^+-\eps$};
\draw[dashed] (2,-2.5)--(-5.5,-2.5);
\draw (2,-2.5) node[anchor=west] {$l_L^-+\eps$};
\draw[line width=.8pt] (-5.5,0.5) .. controls (-5,0.5) and (-4.5,3) .. (-4,3) .. controls (-3.7,3) and (-3.5,-3) .. (-3,-3) .. controls (-2.9,-3) and (-2.1,2.9) .. (-2,2.95) .. controls (-1.9,2.95) and (-1.1,-3.1) .. (-1,-3.1) .. controls (-0.9,-3.1) and (-0.6,3.3) .. (-0.5,3.3) 
.. controls (-.4,3.3) and (-.1,-2.9) .. (0,-2.9) .. controls (.1,-2.9) and (0.4,3.1) .. (0.5,3.1) 
.. controls (.6,3.1) and (.9,-3) .. (1,-3) .. controls (1.1,-3) and (1.5,3.2) .. (1.6,3.2) 
.. controls (1.7,3.2) and (2.4,-0.6) .. (2.5,-0.6);
\fill (.58,2.5) circle (1.5pt);
\draw (0.58,2.5) node[anchor=north west] {$C_1$};
\fill (1,-3) circle (1.5pt);
\draw (1,-3) node[anchor=north west] {$p_1$};
\fill (.06,-2.5) circle (1.5pt);
\draw (.06,-2.5) node[anchor=south west] {$C_2$};
\fill (.5,3.1) circle (1.5pt);
\draw (0.5,3.1) node[anchor=south west] {$p_2$};
\fill (-.4,2.5) circle (1.5pt);
\draw (-.4,2.5) node[anchor=north west] {$C_3$};
\fill (0,-2.9) circle (1.5pt);
\draw (0,-2.9) node[anchor=north west] {$p_3$};
\fill (-.92,-2.5) circle (1.5pt);
\draw (-.92,-2.5) node[anchor=south west] {$C_4$};
\fill (-.5,3.3) circle (1.5pt);
\draw (-0.5,3.3) node[anchor=south west] {$p_4$};
\draw (2.5,0) node[anchor=north west] {$f_\delta$};
\end{tikzpicture}
\caption{A possible situation for Case~I; out of clarity, we have drawn the points $(t,f_\delta(t))$ for $t=C_i$ or $t=p_i$, with $i=1,\,2,\,3,\,4$.}
\label{fig:sect32}
\end{figure}

Let now $\delta<\eta$ be a fixed, small constant, that will be specified in the sequel, and let us consider the smooth function $f_\delta=f\ast \rho_\delta$. As depicted in Figure~\ref{fig:sect32}, we define
\[
C_1 = \max \Big\{ -\eta\leq x\leq 0,\, f_\delta(x)\geq l_L^+-\eps,\, \min\big\{ f_\delta(y),\, x\leq y\leq 0\big\}\leq l_L^-+\eps\Big\}\,.
\]
We can easily notice that, by~(\ref{choiceepsetaI}), $C_1$ belongs to $(-\eta/16,0)$ as soon as $\delta$ has been taken small enough. Indeed, since $\eps<(l_L^+-l_L^-)/4$, there is a set of positive measure in $(-\eta/16,0)$ where $f$ is smaller than $l_L^-+\eps$, and a set of positive measure where $f$ is larger than $l_L^+-\eps$. As a consequence, the same happens for the smooth function $f_\delta$ provided that $\delta$ is sufficiently small, and this proves the claim. We let then $p_1\in \argmin\big\{f_\delta(t),\, t\in [C_1,0]\big\}$, and notice that by construction $f_\delta(p_1)\leq l_L^-+\eps$. Now, we set $C_2 = \max \big\{-\eta\leq x \leq C_1,\, f_\delta(x)\leq l_L^-+\eps \big\}$, which belongs to $(-\eta/8,0)$  if $\delta\ll 1$, and we take $p_2\in \argmax\big\{f_\delta(t),\, t\in [C_2,p_1]\big\}$, which by construction satisfies $f_\delta(p_2)\geq l_L^+-\eps$. Then, we set $C_3=\max \big\{-\eta\leq x \leq C_2,\, f_\delta(x)\geq l_L^+-\eps \big\}$, which again exists if $\delta\ll 1$, and we take $p_3\in\argmin\big\{f_\delta(t),\, t\in [C_3,p_2]\big\}$. We continue this construction in the obvious way, stopping at the first point $p_N$ such that $C_{N+1}$ cannot be defined, because it would be the maximum of an empty set. As observed before, we can have $N$ as large as desired up to have fixed $\delta$ small enough; in particular, just to fix the ideas we will assume that our choice of $\delta$ is so small that $N\geq 10$. Notice that, for every $1\leq i \leq N$, we have that $f_\delta(p_i)\leq l_L^-+\eps$ if $i$ is odd, and $f_\delta(p_i)\geq l_L^+-\eps$ if $i$ is even. In addition, $p_1$ is a minimum point of $f_\delta$ in $[p_2,0]$, $p_2$ is a maximum point of $f_\delta$ in $[p_3,p_1]$, and similarly every $p_j$ is a minimum (resp., maximum) point in $[p_{j+1},p_{j-1}]$ if $3\leq j\leq N-1$ is odd (resp., even). We conclude this case by defining $q_1=-p_2$, which is a maximum point since $f$, and hence $f_\delta$, is symmetric.

\case{II}{If $f(x)=-f(-x),\ h_L>0$ and $l_L^-< \min \{ 0, l_R^-\}$.}
Now, we suppose again that the left jump $h_L$ is strictly positive, but this time we also assume in addition that $f$ is antisymmetric, and moreover that $l_L^-$ is strictly negative and less than $l_R^-$. In the present case, we select $\eps>0$ and $0<\eta<r/4$ so that (still denoting by $\rho$ the symmetric mollifier fixed at the beginning of Section~\ref{sec:second-order}) we have
\begin{equation}\label{choiceepsetaI2}
\begin{split}
\begin{array}{c}
    \bal
    \eps < \frac{(\overline \Lambda-1)h_L}{18\overline \Lambda+6} < \frac{h_L}4\,, \qquad \eps < \frac{l_R^--l_L^-}2 \, \int_{1-\frac{h_L}{4M}}^1 \rho(t)\,dt\,, \eal\\[10pt]
    \bal\L^1(\{l_L^--\eps<f<l_L^++\eps\}\cap (-2\eta,0)) = 2\eta\,.\eal
\end{array}
\end{split}
\end{equation}
Notice that this requirement is stronger than~(\ref{choiceepsetaI}). Now, for some sufficiently small $\delta<\eta$ we define the points $p_i$ exactly as in the first case; notice that in the costruction we did not use the fact that $f$ is symmetric, so everything works now as it did before, the only difference being the different definition of $\eps$ and $\eta$. In this case, we also define the points $q_1=-p_1$ and $q_2=-p_2$; since $f$ is antisymmetric, so is also $f_\delta$ and then the point $q_1$ (resp., $q_2$) is a local maximum since $p_1$ (resp., $p_2$) was a local minimum (resp., maximum).

\par\noindent{\underline{\it Both cases.}}\emph{ Last definitions and selection of the ``good couple''.}\\
Now, we assume to be either in Case~I or in Case~II, and we conclude the present construction by fixing some last quantities, and selecting a ``good couple'' of consecutive points among the $p_i$'s. First of all, calling $D=b-a$, we set
$C(\eta,D,M,g)$ as
\begin{equation}\label{eq:C}
C(\eta,D,M,g) := 20M\left\{|g''|([\eta/2,D])+2\sup_{\eta/2\leq t\leq D}|g'(t)|\right\}\,.
\end{equation}
Concerning this definition, keep in mind that $M=\norm f_{\infty}\geq \norm {f_\delta}_{\infty}$; moreover observe that, since by assumption $g'\in BV_{loc}((0,+\infty))$, then $|g''|([\eta/2,D])$ makes sense because the measure $|g''|$ has a finite value over the compact interval $[\eta/2,D] \comp (0,+\infty)$. In addition, up to possibly decreasing the value of $\delta$, and calling $\bar\gamma = \max\{|p_1-p_2|,\, |q_1-p_1|\}$, we can also assume
\begin{align}\label{eq:delta-choice}
\delta<\frac14\,\min\{|a|,|b|\}\,, &&
\eps \big|g'(\bar\gamma/2)\big|\geq C(\eta,D,M,g)\,.
\end{align}
We take now the smallest index $1\leq j<N-1$ with the property that
\begin{gather}
p_{j+1}-p_{j+2}\geq\frac{p_j-p_{j+1}}2\,,\label{left-cond}\\
p_{j-1}-p_j\geq\frac{p_j-p_{j+1}}2 \qquad \text{if }j>1\,.\label{right-cond}
\end{gather}
Let us show that such an index exists. We begin by checking if $j=1$ works, which is true if $p_2-p_3\geq (p_1-p_2)/2$. If this is the case, then we are already done and $j=1$. Instead, we must have $p_2-p_3<(p_1-p_2)/2<2(p_1-p_2)$, which means that~(\ref{right-cond}) is surely true with $j=2$. Hence, we are done with $j=2$ as soon as also~(\ref{left-cond}) holds, which is true if $p_3-p_4\geq (p_2-p_3)/2$. Again, if this is true we are done and $j=2$, otherwise condition~(\ref{right-cond}) is surely true with $j=3$ and we have to check whether also~(\ref{left-cond}) holds. We proceed recursively, and either we find the claimed existence of the index $j$, or we must arrive to consider an index $j$ such that $p_{j+2}< -\eta/2$. But in this case we would have
\[
p_1+\frac \eta 2< p_1-p_{j+2}  = \sum_{m=1}^{j+1} (p_m-p_{m+1})\leq (p_1-p_2)\sum_{m=1}^{j+1}2^{1-m} \leq 2(p_1-p_2)\,,
\]
which is impossible since by construction $p_1\in (-\eta/16,0)$ and $p_2\in (-\eta/8, 0)$. Hence, the existence of the smallest $1\leq j<N-1$ satisfying~(\ref{left-cond}) and~(\ref{right-cond}) is established. We conclude by selecting the ``good couple''. If $j>1$ or we are in Case~I, then the good couple is simply $\{p_{j+1},p_j\}$. Suppose instead that $j=1$ and we are in Case~II, so that in particular we have defined the point $q_1=-p_1$. Then, the good couple is again $\{p_{j+1},p_j\} = \{p_2,p_1\}$ if the segment $[p_2,p_1]$ is shorter than $[p_1,q_1]$, while otherwise the good couple is $\{p_1,q_1\}$. We point out that the ``good segment'' (that is, the segment whose endpoints are the good couple) is a segment that is not longer than twice the length of its neighbouring segments. Indeed, this is precisely given by~(\ref{left-cond}) and~(\ref{right-cond}) if $j>1$, so we can only consider the case $j=1$. If the good couple is $\{p_1,\, p_2\}$, then the fact that the segment $[p_2,p_1]$ is shorter than the double of $[p_3,p_2]$ is given by~(\ref{left-cond}), and we have to check that it is also shorter than the double of $[p_1,q_1]$. This is true if $f$ is symmetric because then $q_1=-p_2$, so $[p_1,q_1]$ is longer than $[0,q_1]$, which is as long as $[p_2,0]$, hence longer than $[p_2,p_1]$; and it is also true if $f$ is antisymmetric because in this case the good couple is $\{p_1,p_2\}$ only if the segment $[p_2,p_1]$ is shorter than $[p_1,q_1]$. And finally, if the good couple is $\{p_1,q_1\}$, and thus by construction $f$ is necessarily antisymmetric, then by construction we know that $[p_1,q_1]$ is shorter than $[p_2,p_1]$, which is by definition as long as $[q_1,q_2]$.

\section{Proof of Theorem~\mref{continuity}\label{sec:proof}}

In this final section, we are going to prove Theorem~\mref{continuity}. To do so, it is enough to show that there exists an essential limit for $f$ at any point of $(a,b)$; indeed, from the definition itself it follows that $f$ is continuous. More precisely, if for any point $x\in (a,b)$ an essential limit exists, and we call it $g(x)$, then from the definition we immediately get that $g$ is continuous; and then, the functions $g$ and $f$ are readily seen to coincide a.e.\,thanks to the Lebesgue Differentiation Theorem.

Let us then take a generic point of $(a,b)$, and for simplicity of notation let us assume that this point is $0$. We start by showing that if $f$ has left and right (essential) limit at $0$, then they must be equal. Then, a key idea to show that left and right limits actually exist is that we can basically reduce ourselves to consider the special cases when $f$ is even, and when $f$ is odd; indeed, if $f$ has constant potential around $0$, then the same is true for the symmetric and antisymmetric parts $(f(x)\pm f(-x))/2$. To keep the presentation easier, we will consider all the different cases separately. Let us start with the case when left and right limits exist.

\begin{lemma}\label{doubleexistence}
Let $g$ and $f$ be as in Theorem~\mref{continuity}, with $\psi_f$ being a.e. constant in a neighborhood of $0$ (and not necessarily on the whole $[a,b]$). Assume that both the left and the right limit of $f$ at $0$ exist. Then, they must be equal.
\end{lemma}
\begin{proof}
We distinguish the proof in two subcases.
\case{I}{If $\lim_{t\to 0^+} g(t)=+\infty$.}
Let us call for brevity $l_L$ and $l_R$ the left and the right limit of $f$ at $0$, and let us assume by contradiction that they are not equal; without loss of generality, we can think $l_L<l_R$. We fix a positive constant
\begin{equation}\label{howbigeps}
\eps < \frac {l_R-l_L} 5\,.
\end{equation}
Then, there exists $\eta<r/2$ such that $\psi_f$ is a.e. constant in $(-2\eta,2\eta)$ and
\begin{align}\label{boundsf}
f(t) > l_R -\eps \quad \hbox{for a.e. $t\in (0,\eta)$}\,, &&
\big|f(t) - l_L \big|<\eps \quad \hbox{for a.e. $t\in (-\eta,0)$}\,.
\end{align}
For any positive, very small $\delta<\eta$, we know that
\begin{equation}\label{lastnick1}
0 = \psi_f(-\delta) - \psi_f(0) = \int_a^b f(t)\big( g(-\delta-t) - g(-t)\big)\,dt
= \int_a^b f(t)\big( g(t+\delta) - g(t)\big)\,dt\,,
\end{equation}
where the last equality comes since $g$ is symmetric. We evaluate now the last integral in some subintervals. First of all, we have
\[
\int_\eta^b f(t) \big(g(t+\delta)- g(t)\big)\,dt = 
\delta \int_\eta^b f(t) g'(s_{t,\delta})\,dt\,,
\]
where $s_{t,\delta}$ is a suitable point in $(t,t+\delta)$ for any $t\in [\eta,b]$. Recalling that $g\in C^1(\R\setminus\{0\})$ and that $f\in L^\infty$, we immediately deduce that if $\delta$ is sufficiently small then
\begin{equation}\label{lastnick2}
\int_\eta^b f(t) \big(g(t+\delta)- g(t)\big)\,dt \leq C_1 \delta
\end{equation}
for some constant $C_1$, which depends on $a,\,b,\, f,\,g,\,\eps,\,\eta$ but not on $\delta$. In the very same way, up to increasing the constant $C_1$ and decreasing the value of $\delta$, we have
\begin{equation}\label{lastnick3}
\int_a^{-\eta} f(t) \big(g(t+\delta)- g(t)\big)\,dt  \leq C_1 \delta\,.
\end{equation}
Since $g$ is decreasing in $(0,2\eta)$, then $g(t+\delta)-g(t)<0$ for every $t\in (0,\eta)$, and thus by~(\ref{boundsf})
\begin{equation}\label{lastnick4}\begin{split}
\int_0^\eta f(t) \big(g(t&+\delta)- g(t)\big)\,dt \leq \big( l_R-\eps\big) \int_0^\eta g(t+\delta)- g(t)\,dt\\
&= \big( l_R-\eps\big) \bigg(\int_\eta^{\eta+\delta}  g(t)\,dt - \int_0^\delta g(t)\,dt\bigg)
\leq C_2 \delta  - \big( l_R-\eps\big)\int_0^\delta g(t)\,dt\,,
\end{split}\end{equation}
for some constant $C_2$ again depending on $a,\,b,\, f,\,g,\,\eps,\,\eta$.\par

Finally, we have to evaluate the integral in $(-\eta,0)$. The situation is slightly different with respect to the one just considered. Indeed, for every $t\in (0,\eta)$ we have used the fact that $g(t+\delta)-g(t)<0$, which comes since $g$ is decreasing in $(0,2\eta)$. The opposite inequality $g(t+\delta)-g(t)>0$ is not true for every $t\in(-\eta,0)$, but only for $t\in (-\eta,-\delta/2)$; instead, for $t\in (-\delta/2,0)$ we have again $g(t+\delta)-g(t)<0$. As a consequence, by~(\ref{boundsf}) the estimate now reads as
\begin{equation}\label{lastnick5}\begin{split}
\int_{-\eta}^0 f(t) \big(g(&t+\delta) - g(t)\big)\,dt
=\int_{-\eta}^{-\frac\delta 2} f(t) \big(g(t+\delta) - g(t)\big)\,dt
+\int_{-\frac\delta 2}^0 f(t) \big(g(t+\delta) - g(t)\big)\,dt\hspace{-20pt}\\
&\leq (l_L+\eps) \int_{-\eta}^{-\frac\delta 2} g(t+\delta) - g(t)\,dt+(l_L-\eps)\int_{-\frac\delta 2}^0 g(t+\delta) - g(t)\,dt\\
&= (l_L+\eps) \int_{-\eta}^0 g(t+\delta) - g(t) \,dt - 2\eps \int_{-\frac\delta 2}^0 g(t+\delta) - g(t)\,dt\\
&= (l_L+\eps) \bigg(\int_0^\delta g(t) \,dt-\int_{\eta-\delta}^\eta g(t)\,dt\bigg) - 2\eps \bigg(\int_{\frac\delta 2}^\delta g(t)\,dt- \int_0^{\frac\delta 2} g(t)\,dt \bigg)\\
&\leq (l_L+3\eps) \int_0^\delta g(t)\,dt + C_2 \delta \,.
\end{split}\end{equation}
Inserting~(\ref{lastnick2}), (\ref{lastnick3}), (\ref{lastnick4}) and~(\ref{lastnick5}) into~(\ref{lastnick1}), and recalling~(\ref{howbigeps}), we obtain
\[
0\leq 2(C_1+C_2) \delta- \big( l_R-l_L-4\eps\big)\int_0^\delta g(t)\,dt
<2(C_1+C_2) \delta- \eps\int_0^\delta g(t)\,dt\,.
\]
And in turn, this last inequality is impossible as soon as $\delta$ is small enough, because the assumption that $\lim_{t\to 0^+} g(t)=+\infty$ implies that $\int_0^\delta g \gg \delta$.

\case{II}{If $\lim_{t\to 0^+} g(t)<+\infty$.}
Let us now assume that the limit of $g$ at $0$ is finite; hence, $g$ is continuous on the whole $\R$. As a consequence, also recalling that $f\in L^\infty$, the potential $\psi_f$ is differentiable, with
\[
\psi_f'(x)= \int_a^b g'(x-t)f(t)\,dt\,.
\]
Since the value of $\psi_f$ is constant in a neighborhood of $0$, this implies that $\psi_f'$ is constant (and actually $0$) near $0$. We can then argue similarly to how we have done in Case~I, being careful since $g'$ is odd, while $g$ was even. More precisely, we fix again $\eps>0$ and find $\eta$ in such a way that~(\ref{howbigeps}) and~(\ref{boundsf}) are in force together with the fact that $\psi_f'$ is constant in $(-2\eta,2\eta)$. As in~(\ref{lastnick1}), we have
\[
0= \int_a^b f(t)\big( g'(t+\delta) - g'(t)\big)\,dt\,,
\]
and since $g'\in BV_{loc}((0,+\infty))$ as in~(\ref{lastnick2}) and~(\ref{lastnick3}) we get
\[
\bigg|\int_{[a,b]\setminus (-\eta,\eta)} f(t) \big(g'(t+\delta)- g'(t)\big)\,dt \bigg| \leq C_1 \delta\,.
\]
This time, for every $t\in (-\eta,-\delta)\cup (0,\eta)$ we have $g'(t+\delta)-g'(t)>0$, while the opposite inequality is true in $(-\delta,0)$. As a consequence, arguing similarly to how done in~(\ref{lastnick4}) and~(\ref{lastnick5}), this time we get
\[\begin{split}
\int_{-\eta}^\eta f(t) \big(g'(t&+\delta) - g'(t)\big)
 \geq (l_L-\eps)  \int_{-\eta}^{-\delta} g'(t+\delta) - g'(t)\,dt \\
 &\quad + (l_L+\eps)\int_{-\delta}^0 g'(t+\delta) - g'(t)\,dt+ (l_R-\eps)\int_0^\eta g'(t+\delta) - g'(t)\,dt\\
&\geq - C_2\delta + (l_L-\eps) \big(g(0)-g(-\delta)\big)+\Big(2(l_L+\eps)-(l_R-\eps)\Big)\big(g(\delta)-g(0)\big)\\
&= -C_2\delta + \big(g(0)-g(\delta)\big) \Big( l_R-l_L -4\eps\Big)\,.
\end{split}\]
Summarizing, we obtain
\[
0 \geq - (C_1+C_2)\delta + \big(g(0)-g(\delta)\big)\eps\,,
\]
and this is again impossible because the assumptions imply that $g'(t)\to-\infty$ when $t\to 0^+$, and thus
\[
\lim_{\delta\to 0^+} \frac{g(0)-g(\delta)}\delta \to +\infty\,.
\]
\end{proof}

Now, we consider the case when $f$ is even.

\begin{lemma}\label{symmetriccase}
Let $g$ and $f$ be as in Theorem~\mref{continuity}, with $\psi_f$ being a.e. constant in a neighborhood of $0$ (and not necessarily on the whole $[a,b]$). Assume that $a=-b$ and that $f$ is even. Then the limit of $f$ at $0$ exists.
\end{lemma}
\begin{proof}
To get the thesis, it is enough to show that $h_L=0$; indeed, this means that the left limit of $f$ at $0$ exists, but since $f$ is even then in that case also the right limit exists and they coincide. Let us then assume by contradiction that $h_L>0$. Let us define the parameters $\eps,\, \eta$ and $\delta$ as in Section~\ref{subseq:selection-critical-points}, Case~I, having possibily reduced $\eta$ so that $\psi_f$ is a.e. constant in $[-2\eta,2\eta]$. Then, consider the function $f_\delta$, define the corresponding points $p_N < p_{N-1} <\, \cdots\, < p_1< 0 < q_1$ and fix the ``good couple'' $\{p_{j+1},p_j\}$.

We first observe that it is not restrictive to assume that $p_{j+1}$ is a local minimum (or, equivalently, that $j$ is even). Indeed, the proof for the case when $p_{j+1}$ is a local maximum is precisely the same; actually, one can also consider the function $-f$ in place of $f$ and modify the construction of Section~\ref{subseq:selection-critical-points} starting from a point $C_1$ where $f_\delta=l_L^-+\eps$ instead of $l_L^+-\eps$, so that the points $p_i$ remain the same but now they are minima when $i$ is odd and maxima when $i$ is even.

We can now apply Lemma~\ref{lemma:rearrangement-comparison} to the smooth function $f_\delta$ in the interval $[p_{j+1},p_j]$. As discussed in Remark~\ref{rem:change-sign}, this ensures the validity of~(\ref{basecase}), which in our case means that either
\begin{equation}\label{evencase1}
-\int^{p_j}_{p_{j+1}} f_\delta'(t) \, \frac d{dt}\, g(|t-p_{j+1}|)\,dt \geq \frac{|f_\delta(p_j)-f_\delta(p_{j+1})|}2\, \Big(|g'(\gamma)|+|g'(\gamma/2)|\Big)\,,
\end{equation}
or
\begin{equation}\label{evencase2}
\int^{p_j}_{p_{j+1}} f_\delta'(t) \, \frac d{dt}\, g(|t-p_j|)\,dt \geq \frac{|f_\delta(p_j)-f_\delta(p_{j+1})|}2\, \Big(|g'(\gamma)|+|g'(\gamma/2)|\Big)\,,
\end{equation}
where we write for brevity $\gamma=p_j-p_{j+1}$. Notice that, by construction, $\gamma\leq \bar\gamma$, where $\bar\gamma = \max\{|p_1-p_2|,\, |q_1-p_1|\}$ has been defined in Section~\ref{subseq:selection-critical-points}.

We start assuming the validity of~(\ref{evencase1}). Since $p_{j+1}$ is a critical point for $f_\delta$, we have the expression~(\ref{eq:second-derivative}) for $\psi_{f_\delta}''(p_{j+1})$. That quantity must actually be $0$ because by assumption $\psi_f$ is a.e. constant in $(-2\eta,2\eta)$, so  by construction $\psi_{f_\delta}$ is constant in $(-2\eta+\delta,2\eta-\delta)\supseteq (-\eta,\eta)$, which contains $p_{j+1}$. Therefore, we can write
\begin{equation}\label{0=I+J+K}\begin{split}
0&=
\int_a^{p_{j+1}} f_\delta'(t)g'(|p_{j+1}-t|)\,dt
-\int_{p_{j+1}}^{p_j} f_\delta'(t)g'(|p_{j+1}-t|)\,dt
-\int_{p_j}^b f_\delta'(t)g'(|p_{j+1}-t|)\,dt\\
&= I + J + K\,,
\end{split}\end{equation}
where we have called $I,\,J$ and $K$ the three integrals. The validity of~(\ref{evencase1}) allows us immediately to estimate the term $J$, which we will show to be the leading term of the three. In fact, since $f_\delta(p_{j+1})\leq l_L^-+\eps$ and $f_\delta(p_j)\geq l_L^+-\eps$, and since for every $t \in (p_{j+1},p_j)$, we have
\[
\frac d{dt}\, g(|t-p_{j+1}|) = g'(t-p_{j+1})= -|g'(t-p_{j+1})|\,,
\]
then~(\ref{evencase1}) gives
\begin{equation}\label{estiJ}
J \geq \frac{h_L-2\eps}2\, \Big(|g'(\gamma)|+|g'(\gamma/2)|\Big)\,.
\end{equation}
Let us now pass to estimate $I$. To do so, it is convenient to define the set
\[
\Z = \left\{t\in[-\eta,p_{j+1}]: f_\delta(t) < f_\delta(s)\ \forall s\in(t,p_{j+1}]\right\}\,,
\]
\begin{figure}[htbp]  
\begin{tikzpicture}[scale=0.03]
\draw (5, 32) -- (272, 32);
\draw (268, 72) .. controls (256, 152) and (250.398, 130.8814) .. (248.131, 112.695) .. controls (245.864, 94.5085) and (246.932, 79.2543) .. (246.1327, 71.6271) .. controls (245.3333, 64) and (242.6667, 64) .. (240, 76.6667) .. controls (237.3333, 89.3333) and (234.6667, 114.6667) .. (232.6667, 129.3333) .. controls (230.6667, 144) and (229.3333, 148) .. (227.3333, 138.6667) .. controls (225.3333, 129.3333) and (222.6667, 106.6667) .. (218.6667, 89.3333) .. controls (214.6667, 72) and (209.3333, 60) .. (203.3333, 54.6667) .. controls (197.3333, 49.3333) and (190.6667, 50.6667) .. (184.5036, 55.9672) .. controls (178.3406, 61.2678) and (172.6811, 70.5356) .. (170.0145, 87.2023) .. controls (167.3478, 103.869) and (167.6739, 127.9345) .. (163.8369, 137.9672) .. controls (160, 148) and (152, 144) .. (148, 137.3333) .. controls (144, 130.6667) and (144, 121.3333) .. (142.6667, 110) .. controls (141.3333, 98.6667) and (138.6667, 85.3333) .. (132, 76) .. controls (125.3333, 66.6667) and (114.6667, 61.3333) .. (106.6667, 56.6667) .. controls (98.6667, 52) and (93.3333, 48) .. (88, 46.6667) .. controls (82.6667, 45.3333) and (77.3333, 46.6667) .. (72.6667, 53.3333) .. controls (68, 60) and (64, 72) .. (62, 84) .. controls (60, 96) and (60, 108) .. (59.3333, 118.6667) .. controls (58.6667, 129.3333) and (57.3333, 138.6667) .. (51.6841, 142.062) .. controls (46.0348, 145.4573) and (36.0697, 142.9147) .. (4, 88);
\draw[dashed] (244.1859, 66.4183)  -- (211.2604, 66.4189);
\draw[dashed] (194.8322, 51.2731)  -- (194.8322, 32)  -- (194.8322, 32);
\draw[dashed] (194.8322, 51.2731)  -- (97.7754, 51.2726);
\draw[dashed] (84.2467, 46.2031)  -- (84.2467, 32);
\draw[dashed] (97.7754, 51.2726)  -- (97.7754, 32);
\draw[dashed] (211.2604, 66.4189)  -- (211.2604, 32);
\draw[blue, line width = 0.6mm] (84.2467, 32)  -- (97.7754, 32);
\draw[blue, line width = 0.6mm] (194.8322, 32)  -- (211.2604, 32);
\node[text=blue]  at (90, 20) {$\mathcal{Z}$};
\node[text=blue]  at (202, 20) {$\mathcal{Z}$};
\node at (256, 23) {$p_2$};
\node at (240, 23) {$p_3$};
\node at (159, 23) {$p_6$};
\node at (47, 23) {$p_8$};
\node at (20, 23) {$-\eta$};
\fill (253.5, 32) circle (1.5cm);
\fill (244, 32) circle (1.5cm);
\fill (157.3364, 32) circle (1.5cm);
\fill (46.4185, 32) circle (1.5cm);
\fill (19.4185, 32) circle (1.5cm);
\end{tikzpicture}
\caption{A possible situation of a set $\Z$ in the proof of Lemma~\ref{symmetriccase}. Here $j=2$ and the set $\Z$, depicted in blue, is a union of two intervals.}
\label{fig:cancelling}
\end{figure}
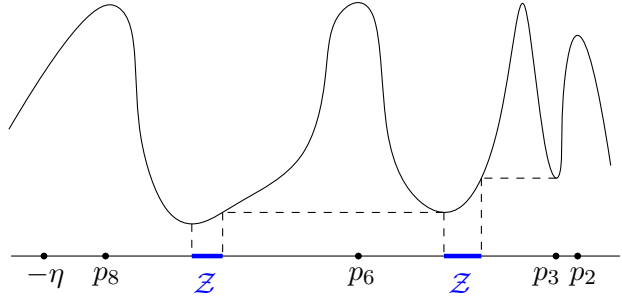
whose meaning appears evident with the aid of Figure~\ref{fig:cancelling}. It is immediate to observe that the restriction of $f_\delta$ to $\Z$ is strictly increasing; moreover, we can write $[-\eta,p_{j+1})\setminus \Z$ as the union of half-open intervals $[a_k,b_k)$, and one has that $f_\delta(a_k)=f_\delta(b_k)$ and $f_\delta'(b_k)=0$. We can then apply Lemma~\ref{lemma:convex-cancellation} to each interval $(a_k,b_k)$, obtaining that
\[
\int_{a_k}^{b_k} f_\delta'(t)g'(p_{j+1}-t)\,dt \geq 0\,,
\]
which implies, also using the change of variables $y=f_\delta(t)$ in $\Z$, that
\begin{equation}\label{cov}
\int_{-\eta}^{p_{j+1}} f_\delta'(t)g'(p_{j+1}-t)\,dt \geq \int_\Z f_\delta'(t)g'(p_{j+1}-t)\,dt
=\int_{f_\delta(\Z)} g'\big(p_{j+1}-f_\delta^{-1}(y)\big)\,dy\,.
\end{equation}
Now, notice that by construction $\Z\subseteq [-\eta, p_{j+2}]$, and by~(\ref{left-cond}) this implies that for every $t\in \Z$ one has $p_{j+1}-t \geq p_{j+1}-p_{j+2}\geq \gamma/2$; moreover, keeping in mind~(\ref{choiceepsetaI}), we know that $f_\delta>l_L^--\eps$ in $(-\eta,0)$ while $f_\delta(p_{j+1})\leq l_L^-+\eps$, and then $\H^1(f_\delta(\Z))\leq 2\eps$; inserting this information in the above estimate gives
\begin{equation}\label{nowdiffer}
\int_{-\eta}^{p_{j+1}} f_\delta'(t)g'(p_{j+1}-t)\,dt \geq 2\eps g'(\gamma/2)=-2\eps |g'(\gamma/2)|\,.
\end{equation}
Moreover, keeping in mind that $g'\in BV_{\rm loc}((0,+\infty))$, we have
\begin{equation}\label{nowequal}
\bigg| \int_a^{-\eta} f_\delta'(t)g'(p_{j+1}-t)\,dt \bigg| \leq \norm {f_\delta}_{\infty} \Big( |g''|([\eta/2,D])+2\sup_{\eta/2\leq t\leq D}|g'(t)|\Big)\leq C(\eta,D,M,g)\,,
\end{equation}
where the constant $C$ has been defined in~(\ref{eq:C}). Since $\eps \big|g'(\bar\gamma/2)\big|\geq C(\eta,D,M,g)$ by~(\ref{eq:delta-choice}), and since $\gamma\leq \bar\gamma$ as noticed before, the last two estimates imply that
\begin{equation}\label{estiI}
I \geq - 3\eps |g'(\gamma/2)|\,.
\end{equation}
In order to estimate the constant $K$ in~(\ref{0=I+J+K}), we can argue in a similar way. More precisely, this time we can call 
\[
\W = \left\{t\in(p_j,\eta]: f_\delta(t) < f_\delta(s)\ \forall s\in(p_j,t)\right\}\,,
\]
notice that $(p_j,\eta)\setminus \W$ is union of half-open intervals $(a_k,b_k]$, apply Lemma~\ref{lemma:convex-cancellation} to obtain that for each such interval
\[
-\int_{a_k}^{b_k} f_\delta'(t)g'(t-p_{j+1})\,dt \geq 0\,,
\]
and argue as in~(\ref{cov}) to get
\[
-\int_{p_j}^\eta f_\delta'(t)g'(|p_{j+1}-t|)\,dt \geq -\int_\W f_\delta'(t)g'(|p_{j+1}-t|)\,dt
=\int_{f_\delta(\W)} g'\big(|p_{j+1}-f_\delta^{-1}(y)|\big)\,dy\,.
\]
The situation now becomes quite different than in the estimate of $I$, and we are going to find a bound which is different from the one of~(\ref{nowdiffer}). In fact, this time we have that $t\geq p_j$ for every $t\in \W$, and then $|p_{j+1}-t|\geq \gamma$, so this time we can pointwise estimate $g'(|p_{j+1}-t|)$ with $g'(\gamma)$, instead of $g'(\gamma/2)$, and this is of course a good news. The bad news, instead, is that it is still true that $f_\delta\geq l_L^--\eps$ in $(0,\eta)$, since $f_\delta$ is a symmetric function, but now the estimate from above coming from~(\ref{choiceepsetaI}) is given by $f_\delta(p_j)\leq l_L^++\eps$, and then the estimate on $\W$ that we get is just $\H^1(f_\delta(W))\leq l_L^+-l_L^-+2\eps=h_L+2\eps$. As a consequence, in place of~(\ref{nowdiffer}) this time we have
\[
-\int_{p_j}^\eta f_\delta'(t)g'(|p_{j+1}-t|)\,dt \geq (h_L+2\eps) g'(\gamma)=-(h_L+2\eps) |g'(\gamma)|\,.
\]
The very same argument as in~(\ref{nowequal}) gives also this time that
\[
\bigg| \int^b_\eta f_\delta'(t)g'(|p_{j+1}-t|)\,dt \bigg| \leq C(\eta,D,M,g)\leq \eps |g'(\gamma/2)|\,,
\]
and then putting everything together we obtain the estimate
\begin{equation}\label{estiK}
K \geq - (h_L+2\eps) |g'(\gamma)|-\eps |g'(\gamma/2)|\,.
\end{equation}
Keeping in mind~(\ref{0=I+J+K}), we now add the estimates~(\ref{estiJ}), (\ref{estiI}) and~(\ref{estiK}) getting
\[\begin{split}
0&=I+J+K \geq \frac{h_L-2\eps}2\, \Big(|g'(\gamma)|+|g'(\gamma/2)|\Big) - 3\eps |g'(\gamma/2)|- (h_L+2\eps) |g'(\gamma)|-\eps |g'(\gamma/2)|\\
&=\frac{-h_L-6\eps}2\, |g'(\gamma)|+\frac{h_L-10\eps}2\, |g'(\gamma/2)| \,.
\end{split}\]
Now, it is time to recall that by~(\ref{goodLambda}) we have $|g'(\gamma/2)| > \overline\Lambda| g'(\gamma)|$, so that the above estimate implies
\[
0\geq \Big(\overline \Lambda(h_L-10\eps) - (h_L+6\eps)\Big)\, \frac{|g'(\gamma)|} 2\,,
\]
which gives the required contradiction thanks to the choice of $\eps$ made in~(\ref{choiceepsetaI}). Summarizing, the contradiction has been found under the assumption  that~(\ref{evencase1}) holds. To conclude the proof, then, we must now work under the assumption that~(\ref{evencase2}) holds. The situation is quite similar to the previous one. In fact, this time we use that $p_j$ is a critical point for $f_\delta$, so we use the expression~(\ref{eq:second-derivative}) with base point $p_j$ instead of $p_{j+1}$, and in place of~(\ref{0=I+J+K}) we get
\[\begin{split}
0&=
\int_a^{p_{j+1}} f_\delta'(t)g'(p_j-t)\,dt
+\int_{p_{j+1}}^{p_j} f_\delta'(t)g'(|p_j-t|)\,dt
-\int_{p_j}^b f_\delta'(t)g'(|p_j-t|)\,dt\\
&= \overline I + \overline J + \overline K\,.
\end{split}\]
The validity of~(\ref{evencase2}) gives, in place of~(\ref{estiJ}), the estimate
\[
-\overline J \geq \frac{h_L-2\eps}2\, \Big(|g'(\gamma)|+|g'(\gamma/2)|\Big)\,.
\]
The estimates for $\overline I$ and $\overline K$, instead, are swapped with respect to those for $I$ and $K$; more precisely, since $p_j$ is a local maximum, the very same arguments as those giving~(\ref{estiI}) and~(\ref{estiK}) imply this time
\begin{align*}
\overline K \leq 3\eps |g'(\gamma/2)|\,, && \overline I \leq  (h_L+2\eps) |g'(\gamma)|-\eps |g'(\gamma/2)|\,.
\end{align*}
Hence, we conclude as in the first case.
\end{proof}

\begin{remark}\label{abitmore}
In the proof of Lemma~\ref{symmetriccase} we have used the symmetry assumption only once, that is, to say that $f_\delta\geq l_L^--\eps$ in $(0,\eta)$ in the estimate of $K$. Indeed, in general we can only say that $f_\delta\geq l_R^--\eps$ in $(0,\eta)$, but the symmetry of $f$ clearly gives $l_L^-=l_R^-$. But then, the proof that $h_L=0$ of Lemma~\ref{symmetriccase} remains true also if $f$ is not symmetric but $l_L^-=l_R^-$, provided the points $p_j$ are defined as in Case~I of Section~\ref{subseq:selection-critical-points}. However, the fact that $h_L=0$ does not imply that $f$ is continuos at $0$ in a general case, while it clearly does so if $f$ is even.
\end{remark}

\begin{lemma}\label{oddcase}
Let $g$ and $f$ be as in Theorem~\mref{continuity}, with $\psi_f$ being a.e. constant in a neighborhood of $0$ (and not necessarily on the whole $[a,b]$). Assume that $a=-b$ and that $f$ is odd. Then the limit of $f$ at $0$ exists.
\end{lemma}
\begin{proof}
We are going to prove that $h_L=0$. Indeed, this guarantees that the left limit of $f$ at $0$ exists, and then also the right limit exists since $f$ is odd. The fact that the two limits coincide, which completes the proof, is then ensured by Lemma~\ref{doubleexistence}.

Let us assume that $h_L>0$, and let us seek for a contradiction. First of all, we notice that since $f$ is odd, then $h_L>0$ implies that $l_L^-\wedge l_R^-<0$. Moreover, as observed in Remark~\ref{abitmore}, if $l_L^-=l_R^-$ then the contradiction is directly given by the same argument as in Lemma~\ref{symmetriccase}, so there is nothing to prove. Considering then the case $l_L^-\neq l_R^-$, we can assume without loss of generality that $l_L^-<l_R^-$, up to replacing $f(x)$ with $f(-x)$. Summarizing, we have to find a contradiction assuming that $h_L>0$ and $l_L^- < \min\{l_R^-,0\}$.\par

In this case, we define the parameters $\eps,\, \eta$ and $\delta$ as in Section~\ref{subseq:selection-critical-points}, Case~II, in such a way that $\psi_f$ is a.e. constant in $[-2\eta,2\eta]$; then, consider the function $f_\delta$, define the corresponding points $p_N < p_{N-1} <\, \cdots\, < p_1< 0 < q_1<q_2$ and fix the ``good couple''. We divide our analysis in three subcases.

\case{I}{The good couple is $\{p_1,\, q_1\}$.}
This is the easiest case to deal with. Indeed, the antisymmetry of $f$ guarantees that
\[
\int_{p_1}^{q_1}f_\delta'(t)|g'(t-p_1)|\,dt = \int_{p_1}^{q_1}f_\delta'(t)|g'(q_1-t)|\,dt\,,
\]
and then~(\ref{eq:cancellation-inequality}) of Lemma~\ref{lemma:rearrangement-comparison} imply that
\[\begin{split}
    -\int_{p_1}^{q_1}f_\delta'(t)  g'(t-p_1)\,dt&=
    \int_{p_1}^{q_1}f_\delta'(t) |g'(t-p_1)|\,dt \geq \frac{f_\delta(q_1)-f_\delta(p_1)}2\, \big(|g'(\gamma)|+|g'(\gamma/2)|\big)\\
    &\geq \frac{l_R^+-l_L^--2\eps}2\, \big(|g'(\gamma)|+|g'(\gamma/2)|\big)
    >\frac{h_L-2\eps}2\, \big(|g'(\gamma)|+|g'(\gamma/2)|\big)\,,
\end{split}\]
where $\gamma=q_1-p_1$ and in the last inequality we have used the fact that $l_R^+>l_L^+$. Notice that this estimate is stronger than the estimate of $J$ in~(\ref{estiJ}). We can then estimate $\int_a^{p_1} f_\delta'(t)g'(|p_1-t|)\,dt$ exactly as in the estimate of $I$ in~(\ref{estiI}). In fact, that estimate uses only that $f\geq l_L^--\eps$ a.e. in $(-\eta,p_1)$, and does not use the symmetry of $f$. Hence, we get
\[
\int_a^{p_1} f_\delta'(t)g'(|p_1-t|)\,dt \geq - 3\eps |g'(\gamma/2)|\,.
\]
And finally, arguing exactly as in the estimate of $\overline I$ in Lemma~\ref{symmetriccase}, we can estimate
\[
\int_{q_1}^b f_\delta'(t)g'(|p_1-t|)\,dt=\int_a^{p_1} f_\delta'(t)g'(|q_1-t|)\,dt\leq (h_L+2\eps) |g'(\gamma)|-\eps |g'(\gamma/2)|\,.
\]
Putting the last three estimates together we get
\[
0=
\int_a^{p_1} f_\delta'(t)g'(|p_1-t|)\,dt
-\int_{p_1}^{q_1} f_\delta'(t)g'(|p_{1}-t|)\,dt
-\int_{q_1}^b f_\delta'(t)g'(|p_{1}-t|)\,dt>0\,,
\]
so the desired contradiction concludes in this case.\par

We can then assume that the good couple is $\{p_j,\, p_{j+1}\}$ for some $j$. This time, we are not allowed to assume that $p_{j+1}$ is a local minimum as we did in Lemma~\ref{symmetriccase}. In fact, to get this assumption one could have to pass from $f$ to $-f$, and this time this is prevented by the assumption that $l_L^-<\min\{l_R^-,0\}$. Let us call again $\gamma=p_j-p_{j+1}$, and let us apply again Lemma~\ref{lemma:rearrangement-comparison} to the smooth function $f_\delta$ in the interval $[p_{j+1},p_j]$. As explained in Remark~\ref{rem:change-sign}, since one between $p_{j+1}$ and $p_j$, call it $p^+$, is an absolute maximum of $f_\delta$ in $[p_{j+1},p_j]$, and the other one, call it $p^-$, is an absolute minimum, we still have~(\ref{basecase}), which in the present case means that either 
\begin{equation}\label{oddcase2}
-\int_{p_{j+1}}^{p_j} f_\delta'(t) \frac d{dt}\, g(|t-p^-|)\,dt \geq \frac{|f_\delta(p_{j+1})-f_\delta(p_j)|}2\, \Big(|g'(\gamma)|+|g'(\gamma/2)|\Big)\,,
\end{equation}
or
\begin{equation}\label{oddcase3}
\int_{p_{j+1}}^{p_j} f_\delta'(t) \frac d{dt}\, g(|t-p^+|)\,dt \geq \frac{|f_\delta(p_{j+1})-f_\delta(p_j)|}2\, \Big(|g'(\gamma)|+|g'(\gamma/2)|\Big)\,.
\end{equation}
We conclude then the proof separately in these two cases.

\case{II}{The good couple is $\{p_{j+1},\,p_j\}$ and~(\ref{oddcase2}) holds.}
This case is also easy, since the situation is close to the one already considered in Lemma~\ref{symmetriccase}. Let us be more precise; we first suppose that~(\ref{oddcase2}) holds, and that $p_{j+1}$ is the minimum point. Then, exactly as in~(\ref{estiJ}) of Lemma~\ref{symmetriccase}, the estimate~(\ref{oddcase2}) gives
\[
-\int_{p_{j+1}}^{p_j} f_\delta'(t)g'(|p_{j+1}-t|)\,dt \geq \frac{h_L-2\eps}2\, \Big(|g'(\gamma)|+|g'(\gamma/2)|\Big)\,.
\]
In addition, exactly as in~(\ref{estiI}) we also obtain
\[
\int_a^{p_{j+1}} f_\delta'(t)g'(|p_{j+1}-t|)\,dt \geq - 3\eps |g'(\gamma/2)|\,,
\]
since the argument used there only considered $f$ in $(a,p_{j+1})\subseteq (a,0)$, and then the fact whether $f$ is symmetric or antisymmetric has no effect. Finally, we need to get
\[
-\int_{p_j}^b f_\delta'(t)g'(|p_{j+1}-t|)\,dt \geq  - (h_L+2\eps) |g'(\gamma)|-\eps |g'(\gamma/2)|\,,
\]
which is analogous to the estimate of the term $K$ in Lemma~\ref{symmetriccase}. In that particular situation, we used the symmetry assumption only to get that $l_L^-=l_R^-$, which yielded that $f_{\delta}\geq l_L^--\eps$ in $(0,\eta)$. However, in the present case, we are assuming that $l_L^-<l_R^-$, and thus we automatically have that $f_{\delta}\geq l_R^--\eps>l_L^--\eps$ in $(0,\eta)$. So the thesis is obtained in this case under the additional assumption that the minimum point is $p_{j+1}$.

Let us now assume that~(\ref{oddcase2}) holds, and that the minimum point is $p_j$. The argument is completely symmetric to the one just performed. Indeed, first of all the validity of~(\ref{oddcase2}) gives
\[
\int_{p_{j+1}}^{p_j} f_\delta'(t)  g'(p_j-t)\,dt \geq \frac{|f_\delta(p_{j+1})-f_\delta(p_j)|}2\, \Big(|g'(\gamma)|+|g'(\gamma/2)|\Big)\,.
\]
Then, in the interval $[p_j,\eta]$ we can estimate
\[
\int_{p_j}^\eta f_\delta'(t)  g'(p_j-t)\,dt \geq -2\eps |g'(\gamma/2)|
\]
using the fact that $f_\delta\geq l_L^--\eps$ in $[p_j,\eta]$ since $l_L^-<l_R^-$. And finally, in the interval $[-\eta,p_{j+1}]$ we have
\[
\int_{-\eta}^{p_j} f_\delta'(t)  g'(p_j-t)\,dt \geq -(h_L+2\eps) |g'(\gamma)|\,,
\]
so that the conclusion is exactly as before.

\case{III}{The good couple is $\{p_{j+1},\,p_j\}$ and~(\ref{oddcase3}) holds.}
This case is the most complicate one, because in the previous case the fact that $l_L^-<l_R^-$ was helping in obtaining the desired estimate, while this time the effect goes against it. To work with this case, we start claiming that
\begin{equation}\label{esti31}
f_\delta(p_i)<l_L^++\eps \quad \forall \, 1\leq i\leq N\,.
\end{equation}
Since $f_\delta(p_i)\leq l_L^-+\eps < l_L^+-\eps$ for all odd $i$, the claim is obvious for those $i$. In addition, since $l_L^--\eps<f<l_L^++\eps$ in $(-2\eta,0)$, the claim is clear also for every odd $i$ such that $p_i<-\delta$. As a consequence, the full validity of~(\ref{esti31}) is established as soon as we prove that
\begin{equation}\label{claim}
p_2<-\delta\,.
\end{equation}
Now, keep in mind that by construction we have $l_L^--\eps<f<l_L^+-\eps$ in $(-2\eta,0)$, and by antisymmetry also $l_L^--\eps<l_R^--\eps<f<l_R^+-\eps$ in $(0,2\eta)$. As a consequence, also recalling that $\delta<\eta$, for every $\delta\big(-1+\frac{h_L}{4M}\big)<t<\eta$ we have
\begin{equation}\label{1533}\begin{split}
    f_\delta(t) &= \int_{-\delta}^\delta f(t+\tau) \rho_\delta(\tau)\,d\tau
    > (l_L^--\eps)\int_{-\delta}^{\delta\big(1-\frac{h_L}{4M}\big)} \rho_\delta(\tau)\,d\tau+(l_R^--\eps)\int_{\delta\big(1-\frac{h_L}{4M}\big)}^\delta \rho_\delta(\tau)\,d\tau\\
    &=l_L^--\eps + \big(l_R^--l_L^-\big) \int_{\delta\big(1-\frac{h_L}{4M}\big)}^\delta \rho_\delta(\tau)\,d\tau
    =l_L^--\eps + \big(l_R^--l_L^-\big) \int_{1-\frac{h_L}{4M}}^1 \rho(\tau)\,d\tau
    >l_L^-+\eps\,,
\end{split}\end{equation}
where the last inequality comes from the second requirement of~(\ref{choiceepsetaI2}). Since $f_\delta(p_1)<l_L^-+\eps$, this means that $p_1\leq -\delta\big(1-\frac{h_L}{4M}\big)$. In addition, since
\[
\norm{f_\delta'}_{\infty}\leq \norm f_\infty\, \norm{\rho_\delta'}_1\leq 2 \norm f_\infty\, \norm{\rho_\delta}_\infty
\leq \frac{2M}\delta\,,
\]
and since $f_\delta(p_2)\geq l_L^+-\eps$, this gives
\[
p_1 - p_2 \geq \frac \delta{2M}\, \big(f_\delta(p_2) - f_\delta(p_1)\big)
> \frac \delta{2M}\, \big(l_L^+-l_L^--2\eps\big)
=\frac \delta{2M}\, \big(h_L-2\eps\big)
>\frac{h_L \delta}{4M}\,.
\]
Putting this estimate together with the one for $p_1$, we have obtained~(\ref{claim}), and as noticed before this implies~(\ref{esti31}).

Let us now seek for a contradiction in this last case. We are assuming that the good couple is $\{p_{j+1},\, p_j\}$ and that~(\ref{oddcase3}) holds. We should subdivide this case in two subcases, namely, whether the maximum point $p^+$ is $p_{j+1}$ or $p_j$. However, exactly as already happened in Case~II, the situation is completely symmetric and the proofs in the two subcases are almost identical. Hence, we only consider what happens if the maximum point is $p_{j+1}$ or, in other words, if $j$ is odd. We try to find a contradiction as already done several times; in fact, the assumption that~(\ref{oddcase3}) holds gives
\[
\int_{p_{j+1}}^{p_j}  f_\delta'(t) g'(t-p_{j+1}) > \frac{h_L-2\eps}2\,  \Big(|g'(\gamma)|+|g'(\gamma/2)|\Big)\,,
\]
while arguing exactly as already done multiple times we obtain the estimates
\begin{align*}
\int_a^{-\eta} f_\delta'(t) g'(p_{j+1}-t)\,dt < \eps |g'(\gamma/2)| \,, && \int_\eta^b f_\delta'(t) g'(t-p_{j+1})\,dt > - \eps |g'(\gamma/2)| \,, \\
\int_{-\eta}^{p_{j+1}}  f_\delta'(t) g'(p_{j+1}-t)\,dt <2\eps |g'(\gamma/2)|\,, &&
\int_{p_j}^{p_1}  f_\delta'(t) g'(t-p_{j+1})\,dt > -(h_L+2\eps) |g'(\gamma)|\,.
\end{align*}
Notice that in order to obtain the third and fourth estimate we have used the fact that $l_L^--\eps<f_\delta<l_L^++\eps$ in $(-\eta,p_1)$, which is true thanks to~(\ref{esti31}). Notice also that the fourth estimate deals with the integral in $[p_j,p_1]$: we would obtain the usual contradiction if we had the estimate in the interval $[p_j,\eta]$, but this is not possible in this case using the techniques exploited before. However, adding the above estimates and using~(\ref{goodLambda}) and the first requirement of~(\ref{choiceepsetaI2}) we obtain
\[\begin{split}
0 &= -\int_a^{p_{j+1}} f_\delta'(t) g'(p_{j+1}-t)\, dt + \int_{p_{j+1}}^b  f_\delta'(t) g'(t-p_{j+1})\,dt\\
&> \frac{h_L-10\eps} 2\, |g'(\gamma/2)| - \frac{h_L+6\eps} 2\, |g'(\gamma)| + \int_{p_1}^\eta  f_\delta'(t) g'(t-p_{j+1})\,dt\\
&> \bigg(\frac{h_L-10\eps} 2- \frac{h_L+6\eps}{2 \overline\Lambda}\bigg)\, |g'(\gamma/2)| + \int_{p_1}^\eta  f_\delta'(t) g'(t-p_{j+1})\,dt\\
&> 4\eps  |g'(\gamma/2)| + \int_{p_1}^\eta  f_\delta'(t) g'(t-p_{j+1})\,dt\,.
\end{split}\]
Thus, the proof is immediately concluded if, by chance,
\begin{equation}\label{bychance}
\int_{p_1}^\eta  f_\delta'(t) g'(t-p_{j+1})\,dt\geq -4\eps  |g'(\gamma/2)| \,.
\end{equation}
If this inequality is not true, then we are not able to find a contradiction using $p_{j+1}$ as ``basepoint''; however, we will find a contradiction using instead $p_1$ as basepoint. Let us be more precise. Applying~(\ref{eq:second-derivative}) to the point $p_1$ and keeping in mind as usual that $\psi_{f_\delta}$ is constant in $(-\eta,\eta)$, we have
\begin{equation}\label{0=0'}
0=-\int_a^{p_1} f_\delta'(t)g'(p_1-t)\,dt+\int_{p_1}^b f_\delta'(t)g'(t-p_1)\,dt\,.
\end{equation}
Now, we can notice that $p_1$ is the minimum point of $f_\delta$ in $[p_1,\eta]$. Indeed, by construction it is the minimum point in $[p_1,0]$; moreover, for every $t\in [0,\eta]$ we have $f_\delta(t)>l_L^-+\eps$ by~(\ref{1533}), while $f_\delta(p_1)<l_L^-+\eps$ by the construction of Section~\ref{subseq:selection-critical-points}. As a consequence, we can apply Lemma~\ref{lemma:concave-cancellation} to the function $F=f_\delta$ with $\alpha=y=p_1$, $\beta=\eta$ and $x=p_{j+1}$ to get that
\[
\int_{p_1}^{\eta} f_\delta'(t)(g'(t-p_{j+1})-g'(t-p_1))\,dt \geq 0\,,
\]
and since~(\ref{bychance}) is false, this gives
\begin{equation}\label{asu0}
\int_{p_1}^{\eta} f_\delta'(t)g'(t-p_1)\,dt <  -4\eps  |g'(\gamma/2)| \,.
\end{equation}
Arguing as usual, we know that
\begin{align}\label{asu1}
\int_a^{-\eta} f_\delta'(t) g'(p_1-t)\,dt > -\eps |g'(\gamma/2)| \,, && \int_\eta^b f_\delta'(t) g'(t-p_1)\,dt < \eps |g'(\gamma/2)| \,.
\end{align}
Finally, arguing once again as in the estimate of $I$ in the proof of Lemma~\ref{symmetriccase}, and keeping in mind that $f_\delta\geq l_L^--\eps$ in $[-\eta,p_1]$ while $f_\delta(p_1)\leq l_L^-+\eps$, we get the estimate
\[
\int_{-\eta}^{p_1} f_\delta'(t)g'(p_1-t)\,dt \geq -2\eps |g'(\gamma')|\,, 
\]
where we have called $\gamma' = p_1-p_2$. By the definition of $j$ made in Section~\ref{subseq:selection-critical-points}, in particular by the fact that $j$ is the smallest value such that~(\ref{left-cond}) and~(\ref{right-cond}) hold, we know that $\gamma'\geq \gamma$, so in particular $|g'(\gamma')|\leq |g'(\gamma)|$ and the above estimate implies
\begin{equation}\label{asu2}
\int_{-\eta}^{p_1} f_\delta'(t)g'(p_1-t)\,dt \geq -2\eps |g'(\gamma)|\,.
\end{equation}
Inserting~(\ref{asu0}), (\ref{asu1}) and~(\ref{asu2}) into~(\ref{0=0'}), we obtain
\[
0=-\int_a^{p_1} f_\delta'(t)g'(p_1-t)\,dt+\int_{p_1}^b f_\delta'(t)g'(t-p_1)\,dt
\leq -2 \eps |g'(\gamma/2)| + 2\eps |g'(\gamma)|<0\,,
\]
which gives the desired contradiction, hence concluding the proof.
\end{proof}

We are finally in position to show our main result, which will be easily obtained by putting together the special cases already considered.
\proofof{Theorem~\mref{continuity}}
As already observed at the beginning of Section~\ref{sec:proof}, a standard argument ensures that $f$ is continuous in $(a,b)$ as soon as it admits an essential limit at every point in that interval.%a standard argument ensures that the continuity of $f$ is guaranteed as soon as it is shown that it admits an essential limit at every point.

Let us then take a generic point $\bar x\in(a,b)$. We define the auxiliary functions $f_S$ and $f_A$ as
\begin{align*}
f_S(x) = f(\bar x+x)+f(\bar x-x)\,, && f_A(x)=f(\bar x+x)-f(\bar x-x)\,.
\end{align*}
Notice that $f_S$ is even and $f_A$ is odd; moreover, since $\psi_f$ is constant in $(a,b)$, then the potential of $f_S$ and that of $f_A$ are constant in $(-s,s)$, where $s=\min\{\bar x-a,b-\bar x\}$. As a consequence, Lemma~\ref{symmetriccase} implies the existence of the essential limit of $f_S$ at $0$, while Lemma~\ref{oddcase} gives the same for $f_A$. The existence of the essential limit of $f$ at $\bar x$ is then obvious.
\end{proof}

%%%%%%%%%%%%%%%%%%%%%%%%%%%%%%%%%%%%%%%%%%%%%%%%%%%%%%%%%%%%%%%%%%%%%%%%%%%%%%%%%%%%
\bigskip
\subsection*{Acknowledgments}
This research was funded in part by the Austrian Science Fund (FWF) [grant DOI \href{https://www.fwf.ac.at/en/research-radar/10.55776/EFP6}{10.55776/EFP6}]. For open access purposes, the authors have applied a CC BY public copyright license to any author-accepted manuscript version arising from this submission.
%%%%%%%%%%%%%%%%%%%%%%%%%%%%%%%%%%%%%%%%%%%%%%%%%%%%%%%%%%%%%%%%%%%%%%%%%%%%%%%%%%%%

\end{document}